\documentclass[12pt]{amsart}
\usepackage{amsmath,amssymb,amsthm,graphicx}

 \newtheorem{thm}{Theorem}[section]

 \newtheorem{conj}[thm]{Conjecture}
 \newcommand{\emb}{\hookrightarrow}

\def\proof{\noindent{\bf Proof. }}
\def\junk#1{}
\begin{document}

\title {Minimal bricks}

\author{Serguei Norine}
\author{Robin Thomas}
\thanks{16 April 2019. Partially supported by NSF grants 0200595 and 0354742. Published in
{\it J.~Combin.\ Theory Ser.~B \bf96} (2006), 505--513. This version fixes an error  kindly
pointed to us by P.~A.\ Fabres, N.~Kothari and M.~H.\ de Carvalho.}
\address{School of Mathematics, Georgia Tech, Atlanta, GA 30332-0160}

\begin{abstract}
A brick is a $3$-connected graph such that the graph obtained from it by
deleting any two distinct vertices has a perfect matching.
A brick is minimal if for every edge $e$ the deletion of $e$ results in
a graph that is not a brick.
We prove a generation theorem for minimal bricks and two corollaries:
(1) for $n \geq 5$, every minimal brick on $2n$ vertices has at most
$5n-7$ edges, and (2) every minimal brick has at least
three vertices of degree three.
\end{abstract}


\maketitle
\section{Introduction}

All the graphs considered in this paper are finite and simple. A \emph{brick}
is a $3$-connected graph such that the graph obtained from it by
deleting any two distinct vertices has a perfect matching. The
importance of bricks stems from the fact that they are building
blocks of the matching decomposition procedure of Kotzig, and
Lov\'asz and Plummer~\cite{LovPlu}. In particular, many matching problems of
interest (such as, for example, computing the dimension of the
linear hull~\cite{EdmLovPul} or lattice~\cite{Lov87} of incidence vectors
of perfect matchings, or characterizing graphs that admit a
``Pfaffian orientation"~\cite{VazYan})
can be reduced to bricks.

In an earlier paper we proved a generation theorem for bricks. The
precise statement requires a large number of definitions, and is given in
Theorem~\ref{extthm} below. Let us describe the result informally first.
Let $G$ be a graph, and let $v_0$ be a vertex of $G$ of degree two
incident with the edges $e_1=v_0v_1$ and $e_2=v_0v_2$. Let $H$ be
obtained from $G$ by contracting both $e_1$ and $e_2$ and deleting
all resulting parallel edges. We say that $H$ was obtained from
$G$ by \emph{bicontracting} or \emph{bicontracting the vertex
$v_0$}, and write $H=G/v_0$.
A subgraph $J$ of a graph $G$ is {\em central} if $G\backslash V(J)$
has a perfect matching. 
We say that a graph $H$ is a {\sl
matching minor} of a graph $G$ if $H$ can be obtained from a central
subgraph of $G$ by repeatedly bicontracting vertices of degree two.
We denote the fact that $H$ is isomorphic to a matching minor of $G$
by writing $H\emb G$.
Our generation theorem of \cite{NorThoBricks} asserts that, 
except for a few well-described exceptions, if $H\emb G$, then
a graph  isomorphic to
$H$ can be obtained from $G$ by repeatedly applying a certain
operation in such a way that all the intermediate graphs are bricks
and no parallel edges are produced.
The  operation is as follows: first delete an
edge, and for every vertex of degree two that results contract
both edges incident with it.
The theorem improves a recent result of 
de~Carvalho, Lucchesi and Murty~\cite{deCLucMurbricks}, but in this
paper we seem to need our result.

We found our theorem useful for generating interesting examples
of bricks and testing various conjectures, but even more useful
was a variant for minimal bricks, which we prove in this paper.
A brick $G$ is {\em minimal} if $G\backslash e$ is not a brick
for every edge $e\in E(G)$. (We use $\backslash$ for deletion.)
The theorem asserts that every minimal brick other than the Petersen
graph can be obtained from $K_4$ or the prism (the complement of a cycle
of length six) by taking ``strict
extensions" in such a way that all the intermediate graphs are
minimal bricks not isomorphic to the Petersen graph. The theorem
is formally stated as Theorem~\ref{minbricks2}.
We postpone the definition of strict extensions until they are needed.

The paper is organized as follows. In the next section we introduce
the results from \cite{NorThoBricks} that we need.
In Section~\ref{sec:generating} we state and prove our generation
theorem for minimal bricks; we deduce it from the more general
Theorem~\ref{extmin}. 
In Section~\ref{sec:edgebound} we prove that, except for
four graphs on at most eight vertices, every minimal brick
on $2n$ vertices has at most $5n - 7$ edges.
Finally, in Section~\ref{sec:3cubicvertices} we prove that every
minimal brick has at least three vertices of degree three.

\section{The tools}

In this section we state the results of~\cite{NorThoBricks} that we need,
but let us start with the following theorem of Lov\'asz~\cite{Loveardec};
see also~\cite[Theorem~5.4.11]{LovPlu}.

\begin{thm}
\label{k4prism}
Every brick has a matching minor isomorphic to $K_4$ or the prism.
\end{thm}

The theorem of de~Carvalho, Lucchesi and Murty~\cite{deCLucMurbricks}
mentioned in the introduction uses $K_4$ and the prism as the starting
graphs of their generation procedure. We use a more restricted
set of operations, and the price we pay for that is that the starting set has
to be expanded. We now introduce the relevant classes of graphs.

Let $C_1$ and $C_2$ be two vertex-disjoint cycles of length
$n\ge3$ with vertex-sets $\{u_1,u_2$,$\ldots$, $u_n\}$ and
$\{v_1,v_2$,$\ldots,v_n\}$ (in order), respectively, and let $G_1$
be the graph obtained from the union of $C_1$ and $C_2$ by adding
an edge joining $u_i$ and $v_i$  for each $i=1,2,\ldots,n$. We say
that $G_1$ is a \emph{planar ladder}. Let $G_2$ be the graph
consisting of a cycle $C$ with vertex-set
$\{u_1,u_2$,$\ldots$,$u_{2n}\}$ (in order), where $n \geq 2$ is an
integer, and $n$ edges with ends $u_i$ and $u_{n+i}$ for
$i=1,2$,$\ldots,n$. We say that $G_2$ is a \emph{M\"obius ladder}.
A \emph{ladder} is a planar ladder or a M\"obius ladder. Let $G_1$
be a planar ladder as above on at least six
 vertices, and let $G_3$ be obtained from $G_1$ by deleting
the edge $u_1u_2$ and contracting the edges $u_1v_1$ and $u_2v_2$.
We say that $G_3$ is a \emph{staircase}. Let $t\ge2$ be an
integer, and let $P$ be a path with vertices
$v_1,v_2,\ldots,v_{t}$ in order. Let $G_4$ be obtained from $P$ by
adding two distinct vertices $x,y$ and edges $xv_i$ and $yv_j$ for
$i=1,t$ and all even $i\in\{1,2,\ldots,t\}$ and $j=1,t$ and all
odd  $j\in\{1,2,\ldots,t\}$. Let $G_5$ be obtained from $G_4$ by
adding the edge $xy$. We say that $G_5$ is an \emph{upper
prismoid}, and if $t\ge4$, then we say that $G_4$ is a \emph{lower
prismoid}. A \emph{prismoid} is a lower prismoid or an upper
prismoid.

We need the following  strengthening of Theorem~\ref{k4prism}, proved
in~\cite[Theorem (1.8)]{NorThoBricks}.

\begin{thm}\label{minbricks}Let $G$ be a brick
not isomorphic to $K_4$, the prism or the Petersen graph. Then $G$
has a matching minor isomorphic to one of the following seven
graphs: the graph obtained from the prism by adding an edge, the
lower prismoid on eight vertices, the staircase on eight vertices,
the staircase on ten vertices, the planar ladder on ten vertices,
the wheel on six vertices, and the M\"obius ladder on eight
vertices.
\end{thm}

In the introduction we described our generation theorem by means
of operations that reduce the larger graph $G$ to its matching
minor $H$. This version is easier to describe concisely, but
for both the proof and the applications it is better to proceed
the other way, namely to describe how to obtain $G$ from $H$. Thus we
reverse the process now and proceed in the other direction.
Here are the relevant definitions.

Let $H,G,v_0,v_1,v_2,e_1,e_2$ be as in the definition of
bicontraction. Assume that $v_1$, $v_2$ are not adjacent, that they
both have degree at
least three and that they have no common neighbors except $v_0$;
then no parallel edges are produced during the contraction of
$e_1$ and $e_2$. Let $v$ be the new vertex that resulted from the
contraction. 
We say that $G$ was obtained from $H$ by \emph{bisplitting
the vertex $v$}. We call $v_0$ the \emph{new inner vertex} and
$v_1$ and $v_2$ the \emph{new outer vertices}. Let $H$ be a graph.
We wish to define a new graph $H''$ and two vertices of $H''$.
Either $H''=H$ and $u,v$ are two nonadjacent vertices of $H$, or
$H''$ is obtained from $H$ by bisplitting a vertex, $u$ is the new
inner vertex of $H''$ and $v\in V(H'')$ is not adjacent to $u$, or
$H''$ is obtained by bisplitting a vertex of a graph obtained from
$H$ by bisplitting a vertex, and $u$ and $v$ are the two new inner
vertices of $H''$. 
Finally, let $H'$ be obtained from $H''$ by adding an edge with
ends $u,v$. We say that $H'$
is a \emph{linear extension} of $H$.

Since in the next theorem the graph $H$ need not be a brick
we need two more exceptional classes of graphs.
Let $C$ be an even cycle with
vertex-set $v_1,v_2, \ldots, v_{2t}$ in order, where $t \geq 2$
is an integer and let $G_6$ be obtained from $C$ by adding
vertices $v_{2t+1}$ and $v_{2t+2}$ and edges joining $v_{2t+1}$ to
the vertices of $C$ with odd indices and $v_{2t+2}$ to the
vertices of $C$ with even indices. Let $G_7$ be obtained from
$G_6$ by adding an edge $v_{2t+1}v_{2t+2}$. We say that $G_7$ is an
{\em upper biwheel}, and if $t\ge3$ we say that $G_6$ is a {\em lower biwheel}.
A {\em biwheel} is a lower biwheel or an upper biwheel.
Please note that biwheels are bipartite, and therefore are not bricks.

We are now ready to state a version of our generation
theorem~\cite[Theorem (1.10)]{NorThoBricks}. The version mentioned
in the introduction follows easily, because a linear extension
of a brick is a brick.

\begin{thm}\label{extthm}Let $G$ be a brick other than the Petersen graph, and
let $H$ be a $3$-connected matching minor of $G$. Assume that if
$H$ is a planar ladder, then
there is no strictly larger planar ladder $L$ with $H\emb L\emb G$,
and similarly for M\"obius ladders, wheels,
lower biwheels, upper biwheels, staircases, lower prismoids and upper prismoids.
If $H$ is not isomorphic to $G$, then
some matching minor of $G$ is isomorphic to a linear extension of $H$.
\end{thm}

\section{Generation Theorem for Minimal Bricks}
\label{sec:generating}

In this section we prove a generation theorem for minimal bricks,
Theorem~\ref{minbricks2} below. We derive it from the more general
Theorem~\ref{extmin}.

If $H$ is a graph, and $u,v\in V(H)$ are distinct nonadjacent
vertices, then
$H+(u,v)$ or $H+uv$ denotes the graph obtained from $H$ by adding
an edge with ends $u$ and $v$. If $u$ and $v$ are adjacent or equal then
$H+uv = H$. Now let $u,v\in V(H)$ be adjacent. By
\emph{bisubdividing} the edge $uv$ we mean replacing the edge by a
path of length three, say a path with vertices $u,x,y,v$, in
order. Let $H'$ be obtained from $H$ by this operation. We say
that $x,y$ (in that order) are the \emph{new vertices}. Thus $y,x$
are the new vertices resulting from subdividing the edge $vu$ (we
are conveniently exploiting the notational asymmetry for edges).
Now if $w\in V(H)-\{u\}$, then by $H+(w,uv)$ we mean the graph
$H'+(w,x)$. Notice that the graphs $H+(w,uv)$ and $H+(w,vu)$ are
different. 


Let $H$ be a graph, let $u,v\in V(H)$ be distinct, and let $H'$ be obtained
from $H+uv$ by bisubdividing $uv$, where the new vertices are $x,y$.
Let $x'\in V(H)-\{u\}$
and $y'\in V(H)-\{v\}$ be not necessarily distinct vertices such
that not both belong to $\{u,v\}$. In those circumstances we say
that $H'+(x,x')+(y,y')$ is a \emph{quasiquadratic extension of $H$}. 
We say that it is a \emph{quadratic extension of $H$} if $u$ and $v$
are not adjacent in $H$. (Recall our convention that
if $u$ and $v$ are adjacent in $H$, then $H + uv = H$.)
We say that $uv$ is \emph{the base} of this quasiquadratic extension. 

Now let $u,v,H',x,y$ be as above, and let 
$a,b\in V(H)$  be  not necessarily distinct vertices such that $\{u,v\} \neq \{a,b\}$,
and if $a=b$ then $a\not \in\{u,v\}$.
If $a\ne b$,  then let $H''$ be obtained from $H'+ ab$ by bisubdividing $ab$, and 
let $x',y'$
be the new vertices. If $a=b$, then let $H''$ be obtained from $H'$ by adding new 
vertices $x',y'$ and edges $ax'$, $x'y'$ and $y'a$.
Then the graph $H''+(x,x')+(y,y')$ is called
a \emph{quasiquartic extension of $H$}. 
It is a \emph{quartic extension of $H$} if $uv,ab\in E(H)$.
We say that $uv, ab$ are
\emph{the bases} of the quasiquartic extension.
%
Quadratic and quartic extensions were used in the proof of 
Theorem~\ref{extthm} in~\cite{NorThoBricks}; quasiquadratic
and quasiquartic extensions are new.

\begin{figure}
\centering
\includegraphics{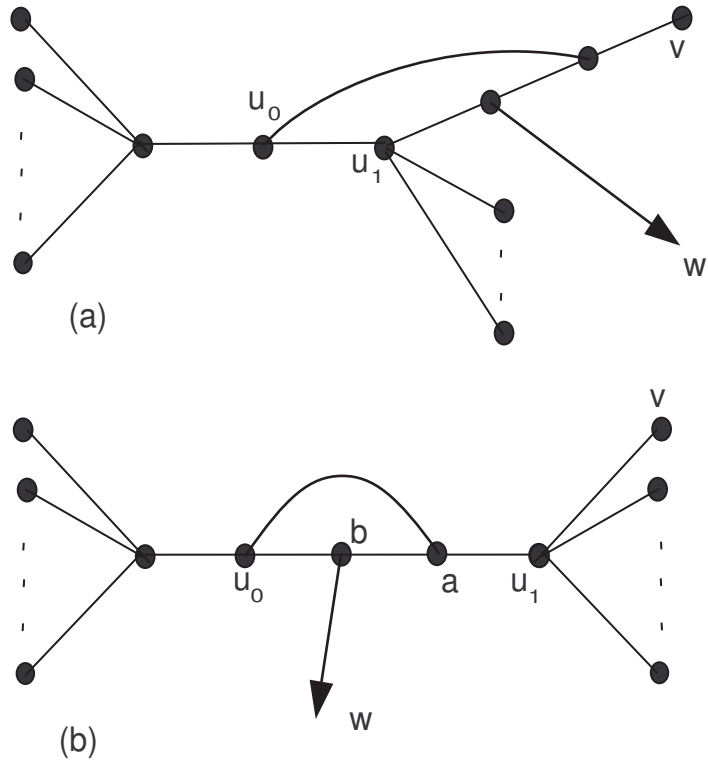}
\caption{(a) Bilinear extension, (b) Pseudolinear extension}
\label{bilinear}
\end{figure}

We need to define two new types of extension.
We say that a linear 
extension $H'$ of a
graph $H$ is \emph{strict} if $|V(H')| > |V(H)|$. Let $u,v,w$ be
pairwise distinct vertices of $H$, let $H'$ be obtained from $H$
by bisplitting $u$, and let $u_0$ be the new inner vertex and $u_1$
a new outer vertex. If $u_1v \in E(H')$ and $vw \not \in E(H)$
then the graph $H' + (u_0,vu_1) + (y, w)$, where $x,y$ are
the new vertices of $H' + (u_0,vu_1)$, is called a
\emph{bilinear} extension of $H$. If $uw \not \in E(H)$ then the
graph $H' +(u_0,u_1u_0)+(b, w)$, where $a,b$ are the new vertices
of $H' +(u_0,u_1u_0)$, is called a
\emph{pseudolinear} extension of $H$. See Figure~\ref{bilinear}.

Finally, we say that $H'$ is a
\emph{strict} extension of $H$ if $H'$ is a quasiquadratic,
quasiquartic, bilinear, pseudolinear or strict linear extension of~$H$.
It is not hard to see that a strict extension of a brick is a  brick.

\begin{thm}\label{extmin}Let $G$ be a brick other than the Petersen graph, and
let $H$ be a $3$-connected matching minor of $G$ such that
$|V(H)| < |V(G)|$. Then some matching minor of $G$ is isomorphic
to a strict extension of $H$.
\end{thm}

\proof Let a graph $H'$ be chosen so that $H$ is a spanning
subgraph of $H'$, $H' \emb G$ and $|E(H')|$ is maximal.

Suppose first that $H'$ is a planar ladder and there exists a
planar ladder $L$ with $H'\emb L\emb G$ and $|V(L)| > |V(H')|$.
Then clearly $H'=H$, and if we choose $L$ with $|V(L)|$ minimum, 
then $L$ is a quartic extension of $H$ and
therefore the theorem holds. Therefore we can assume that if $H'$
is a planar ladder, then there is no strictly larger planar ladder
$L$ with $H\emb L\emb G$, and similarly for M\"obius ladders,
wheels, lower biwheels, upper biwheels, staircases, lower
prismoids and upper prismoids. By Theorem~\ref{extthm} and the
choice of $H'$ there exists a strict linear extension $K$ of $H'$
such that $K \emb G$. We denote $E(H')-E(H)$ by $E'$ and break the
analysis into cases depending on the type of strict linear
extension.

Suppose first that $K = K'+uv$, where $K'$ is obtained from $H'$
by bisplitting a vertex, $v$ is the new inner vertex of $K'$ and
$u \in V(H')$. Let $v_1$ and $v_2$ be the new outer vertices. We
have $E(H') \subseteq E(K')$, in the natural way. For $i=1,2$ let
$d_i$ be the number of edges of $E(H)$ that are incident with
$v_i$ in $K'$ (or $K$). We assume without loss of generality that
$d_1 \geq d_2$. Note that $d_1 + d_2 \geq 3$, because $v$ has
degree at least three in $H$.

If $d_2 \geq 2$ then $K\setminus E'$ is a strict linear extension
of $H$. If $d_2=1$ let $f \in E'$ be an edge incident with $v_2$;
then $K \setminus (E' - \{f\})$ is a quadratic extension of $H$.
Finally, if $d_2=0$ and $f_1, f_2 \in E'$ are incident with $v_2$
then $K \setminus (E' - \{f_1,f_2\})$ is a quasiquadratic
extension of $H$.

Now suppose $K = K'+ u_1u_2$, where $K'$ is obtained by
bisplitting a vertex of a graph obtained from $H'$ by bisplitting
a vertex, and $u_1$ and $u_2$ are the two new inner vertices of
$K'$. Let $v_1,v_2$ and $v_3,v_4$, respectively, be the
corresponding new outer vertices. Let $d_1,d_2,d_3$ and $d_4$ be
defined analogously as above. We start by assuming that
$v_1,v_2,v_3$ and $v_4$ are pairwise distinct and without loss of
generality assume $d_1 \geq d_2$, $d_3 \geq d_4 \geq d_2$.

If $d_2 \geq 2$ then $K\setminus E'$ is a strict linear extension
of $H$. If $d_2 = 1, d_4 \geq 2$ then $K\setminus E'/v_2$ is
isomorphic to a strict linear extension of $H$ unless the edge of
$H$ incident with $v_2$ is incident also with one of the vertices
$v_3$ and $v_4$. In this case $K\setminus (E' - \{f\})$ is a
bilinear extension of $H$, for every $f \in E'$ incident with
$v_2$. If $d_2 = d_4 = 1$ for $i \in \{1,2\}$ let $e_i$ denote the
unique edge in $E(H)$ incident with $v_{2i}$ and let $f_i$ denote
some edge in $E'$ incident with $v_{2i}$. If $e_1 = e_2$ then
$K\setminus (E' - \{f_1,f_2\})$ is a quasiquartic extension of
$H$. (If $f_1$ is adjacent to $f_2$, then we need the provision of $a=b$ 
in the definition of quasiquartic extension.)
Otherwise, without loss of generality we assume that $e_2$ is
not incident with $v_1$ and deduce that $K\setminus (E' -
\{f_1\})/v_4$ is a quadratic extension of $H$ with base $e_1$.

It remains to consider the subcase when $d_2 = 0$. Let $f,f' \in
E'$ be incident with $v_2$ such that $f$ has no end in
$\{v_3,v_4\}$. If $d_4 \geq 2$ then $K\setminus (E' -
\{f\})\setminus u_1v_1 /u_1$ is a strict linear extension of $H$.
If $d_4 = 1$ let $e$ denote the unique edge in $E(H)$ incident
with $v_4$. If $e$ is not incident with $v_1$ then $K\setminus (E'
- \{f,f'\})/v_4$ is a quasiquadratic extension of $H$ if $f'$ is
not incident with $v_4$ and $K\setminus (E' - \{f,f'\})$ is a
quasiquartic extension of $H$ if $f'$ is incident with $v_4$. If
on the other hand $e$ is incident with $v_1$ then $K\setminus (E'
- \{f,f''\})\setminus u_1v_1 /u_1$ is a quadratic extension of
$H$, where $f''$ is any edge in $E'$ incident with $v_4$. Finally,
if $d_4 = 0$ let $f^* \in E'$ be incident with $v_4$ and have no
end in $\{v_1,v_2\}$. Then $K\setminus (E' -
\{f,f',f^*\})\setminus u_2v_3/u_2$ is a quasiquadratic extension
of $H$. This completes the case when $v_1,v_2,v_3$ and $v_4$ are
pairwise distinct.

We now assume without loss of generality that $v_1 = v_4$. Then
$v_1,v_2$ and $v_3$ are pairwise distinct and we assume $d_2 \geq
d_3$, again without loss of generality. Suppose first $d_1 = 0$.
If $d_3 \geq 2$ then $K\setminus (E' - \{g\})$ is a pseudolinear
extension of $H$, where $g \in E'$ is incident with $v_1$;  if
$d_3 = 1$ then $K\setminus (E' - \{g\})/v_3$ is a quadratic
extension of $H$ and if $d_3=0$ then $K\setminus (E' -
\{f,g\})/v_3$ is a quasiquadratic extension of $H$, where $f$ is
an edge in $E'$ incident with $v_3$ and not adjacent to $g$.
Therefore we may assume $d_1 \geq 1$. If $d_2 \geq 2$ and $d_3
\geq 1$ then $K \setminus E'$ or $K\setminus E'/v_3$ is a strict
linear extension of $H$. If $d_2 \geq 2$ and $d_3 = 0$ then
$K\setminus (E' \setminus f)/v_3$ is a quadratic extension of $H$,
where $f$ is as above. If, finally, $d_2 \leq 1$ then let $E''$ be
obtained from $E'$ by deleting $2 - d_2$ edges of $E'$ incident
with $v_2$ and $1-d_3$ edges incident with $v_3$; in that case
$K\setminus E''\setminus v_1u_2/u_2$ is a quasiquadratic extension of $H$.

This completes the case analysis.~\qed

\vskip 10pt

Theorem~\ref{extmin} implies the following generation theorem for
minimal bricks.

\begin{thm}\label{minbricks2} Let $G$ be a minimal brick other than the Petersen
graph. Then $G$ can be obtained from $K_4$ or the prism by taking
strict extensions, in such a way that all the intermediate graphs
are minimal bricks not isomorphic to the Petersen graph.
\end{thm}

\proof Suppose the statement of the theorem is false and let $G$
be a counterexample with $|V(G)|$ minimum.

By Theorem~\ref{k4prism} we may choose a minimal brick $H \emb G$ 
such that $H$ can be obtained from $K_4$ or the prism by taking 
strict extensions and, subject to that,
$|V(H)|$ is maximum. 
If $|V(H)|=|V(G)|$ then $H$ is isomorphic
to $G$ by the minimality of $G$. If, on the other hand, $|V(H)| <
|V(G)|$, then by Theorem~\ref{extmin} there exists a strict
extension $H' \emb G$ of $H$. Let $H'' \emb H'$ be a minimal brick
with $|V(H'')|=|V(H')|$; then $H'' \emb G$. It follows that $H''$
is not isomorphic to $G$, for otherwise so is $H'$, contrary to
our assumption that $G$ is a counterexample to the theorem. By the
minimality of $G$ the graph $H''$ can be obtained from $K_4$ or
the prism by taking strict extensions, contrary to the choice of
$H$.~\qed

\vskip 10pt

Note that there exist bricks obtained from $K_4$ or the prism by a
sequence of strict extensions, that are not minimal. A simple
example follows.

Let $G$ be the prism, $V(G)=\{v_1,v_2,v_3,u_1,u_2,u_3\}$, the
vertices $v_1$,$v_2$,$v_3$ are pairwise adjacent and so are the
vertices $u_1,u_2,u_3$, and $u_i$ is adjacent to $v_i$ for $i \in
\{1,2,3\}$. Let $G'=G+u_1v_2$ and let $G''=G' + (u_2,u_1v_2) +
v_1y$, where $x,y$ are the new vertices of $G' + (u_2,u_1v_2)$. 
Then $G''$ is a quasiquadratic extension of $G$ and
$G'' \setminus u_1v_1$ is a brick, which can be obtained from a
prism by a quadratic extension or a sequence of two linear
extensions.

\section{Edge Bound for Minimal Bricks}
\label{sec:edgebound}

The following theorem is~\cite[Corollary~5.4.16]{LovPlu}.

\begin{thm} If $G$ is a minimal bicritical graph with $n \geq
6$ vertices, then $|E(G)| \leq 5(n-2)/2$.
\end{thm}

We use Theorem~\ref{extmin} to prove a similar bound for minimal bricks.

\begin{thm}\label{edgebound} Let $G$ be a
minimal brick on $2n$ vertices. Then $|E(G)| \leq 5n - 7$, unless
$G$ is the prism or the wheel on four, six or eight vertices.
\end{thm}

\proof The theorem holds for the Petersen graph, so from now on we
assume that $G$ is not the Petersen graph, the prism or the wheel
on six or eight vertices. Denote the last three graphs by $R_6$,
$W_6$ and $W_8$, respectively.

Note that a strict linear extension increases the number of
vertices in a graph by $2$ or $4$ and the number of edges by $3$
or $5$, respectively. Similarly, a quasiquadratic extension
increases the number of vertices by $2$ and the number of edges by
at most $5$, while quasiquartic, bilinear and pseudolinear
extensions increase the number of vertices by $4$ and the number
of edges by at most $8$.

We say that a brick $H$ is \emph{sparse} if $|E(H)| \leq
\frac{5}{2}|V(H)| - 7$ and we say that $H$ is {\em dense} otherwise. We
claim that any minimal brick that contains a sparse matching minor is
sparse. Suppose $G_1$ and $G_2$ are bricks, $G_1 \emb G_2$,  $G_1$
is sparse and $G_2$ is minimal. Let a sparse brick $H \emb G_2$
be chosen with $|V(H)|$ maximum. From Theorem~\ref{extmin} we
deduce that either $|V(H)|=|V(G_2)|$ or some strict extension $H'$
of $H$ is a matching minor of $G_2$. In the latter case, by the
calculations above, $H'$ is sparse in contradiction with the
choice of $H$. Therefore $|V(H)|=|V(G_2)|$ and $G_2$ is isomorphic
to $H$ by the minimality of $G_2$. The claim follows.

Suppose $G$ is dense. By Theorem~\ref{minbricks} $G$ has a
matching minor isomorphic to one of the seven graph mentioned therein,
and hence $G$ has a matching minor isomorphic to one of the
following four graphs:
$R_6$, $W_6$, the staircase on eight vertices, and the M\"obius
ladder on eight vertices. Among these graphs only two are dense:
$R_6$ and $W_6$.

\begin{figure}
\centering
\includegraphics{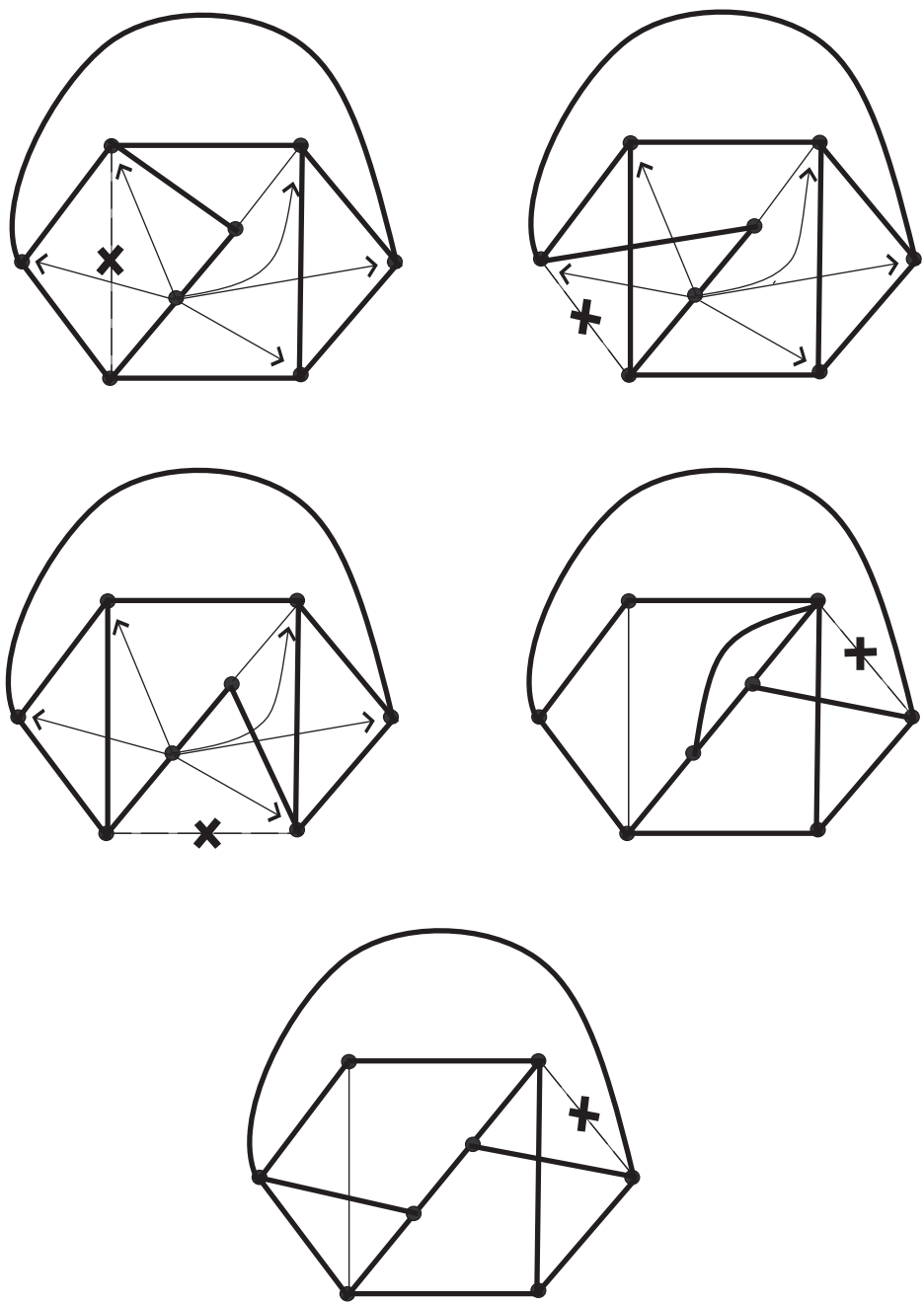}
\caption{Quasiquadratic extensions of the prism with $14$ edges}
\label{QUASI}
\end{figure}

Assume first that $G$ contains $R_6$ as a matching minor. By
Theorem~\ref{extmin} there exists a strict extension $H$ of the
prism such that $H \emb G$. By the calculations above $H$ is
sparse, unless $H$ is a quadratic extension of $R_6+uv$ with
base $uv$, where $uv \not \in E(R_6)$. We will show that there
exists $e \in E(H)$ such that $H\setminus e$ is a brick. Note that
$H\setminus e$ is sparse. Therefore it follows that any minimal
brick containing the prism as a matching minor and not equal to it
is sparse. We prove the existence of $e$ by listing all possible
quasiquadratic extensions of $R_6$ with $14$ edges in
Figure~\ref{QUASI}. An edge $e$ that satisfies the conditions
above is indicated by a cross. A spanning bisubdivision or bisplit
of $R_6$ or $W_6$ in $H\setminus e$ is indicated by bold lines and
allows the reader to easily verify that the claim holds in each of
the cases.

Therefore we may assume that $G$ contains $W_6$ as a matching
minor and does not contain $R_6$. By Theorem~\ref{extthm} $G$ is
a wheel or $G$ contains a linear extension of a wheel as a matching
minor. All the wheels on at least ten vertices and all strict
linear extensions of $W_6$ and $W_8$ are sparse and therefore $G$ must
contain a graph obtained from $W_6$ or $W_8$ by an edge addition. 
Every graph obtained from $W_8$ by adding an edge has a matching
minor isomorphic to a graph obtained from $W_6$ by adding an edge.
The latter
graph is unique up to isomorphism and contains $R_6$ as a spanning
subgraph, in contradiction with our assumptions. ~\qed

\vskip 10pt

The bound given in Theorem~\ref{edgebound} is tight for every $n
\geq 4$. An example of a minimal brick $G_n$ on $2n + 4$ vertices
with $5n + 3$ edges for $n \geq 2$ follows. Let
$V(G_n)=\{x,y,z,t,v_1,u_1,v_2,u_2, \ldots, v_n,u_n\}$. For every
$i \in \{1,2, \ldots, n\}$ let $xt,yt,zt,xu_i$,$yu_i,yv_i$,$zv_i$
and $u_iv_i$ be the edges of $G_n$. Then for every $e \in E(G_n)$
the graph $G_n\setminus e$ contains a vertex of degree two, and hence is
not a brick. It remains to show that $G_n$ is a brick for every
$n$. Note that $G_k$ is a quasiquadratic extension of $G_{k-1}$
for every $k>2$. Therefore it suffices to show that $G_2$ is a
brick. The graph $G_2\setminus u_1y \setminus v_1y$ is isomorphic
to the prism with one of its edges bisubdivided and consequently
$G_2$ can be obtained from the prism by a quadratic extension.

\section{Three Cubic Vertices}
\label{sec:3cubicvertices}


De Carvalho, Lucchesi and Murty~\cite{deCLucMurbricks} proved
that every minimal brick has a vertex of degree three.
According to them (private communication) it had been conjectured
by Lov\'asz. We prove the following strengthening.

\begin{thm}\label{threevertices} Every
minimal brick has at least three vertices of degree three.
\end{thm}

\proof Let a minimal brick $G$ that has at most two vertices of
degree three be chosen with $|V(G)|$ minimal. By
Theorem~\ref{minbricks2} there exists a minimal brick $H \emb G$
with at least three vertices of degree three, such that $G$ is
isomorphic to a strict extension of $H$.

Note that if a strict linear extension is used to obtain $G$ from
$H$ then the degree of at most one vertex of $H$ increased and at
least one vertex in $V(G) - V(H)$ has degree three. If a
quasiquartic, bilinear or pseudolinear extension is used to obtain
$G$ then $V(G) - V(H)$ contains at least three vertices of degree
three. Therefore $G$ is isomorphic to a quasiquadratic extension
of $H$ that is not quadratic.

We assume without loss of generality that $V(G) -
V(H)=\{u_1,u_2\}$ and there exist $v_1, v_2, v_3, v_4 \in V(H)$
such that $E(G) - E(H)= \{u_1v_1,$ $u_1v_2$, $u_2v_3,$
$u_2v_4,u_1u_2\}$, at least three of the vertices $v_1, v_2, v_3,
v_4$ are distinct, $v_1 \neq v_2$ and $v_3 \neq v_4$. Note that
the vertices of degree three in $H$ must form a subset of
$\{v_1,v_2,v_3,v_4\}$ and that $v_1v_3,v_2v_3,v_2v_4,v_1v_4 \not
\in E(H)$, for the deletion of such an edge from $G$ results in a quadratic
extension of $H$, contrary to the fact that $G$ is a minimal
brick.

Since $H$ is a brick, it is not a biwheel.
By Theorem~\ref{extthm} either $H$ is a ladder, wheel, staircase
or prismoid or $H$ is a linear extension of a brick. 
If $H$ is a ladder, wheel,
staircase or prismoid distinct from $K_4$ then $H$ has at least
$5$ vertices of degree three, and consequently $G$ has at least
three vertices of degree three. If $H=K_4$ then $G$ is not
minimal, by an observation in the previous paragraph.

Therefore, $H$ is a linear extension of a brick, and hence
 there exists $e \in E(H)$ such that $H\setminus e$
becomes a brick after possible bicontractions of vertices of
degree two in such a way that no parallel edges are created by these
bicontractions. Note that $H$ is minimal and therefore at least
one end of $e$ is a vertex of degree three in $H$. Assume first
that exactly one end of $e$ has degree three in $H$. Without loss
of generality this end is $v_1$. The graph $G \setminus e$ is a
brick, because it can be obtained by a linear extension (first
bisplit to produce $H \backslash e$, then add the edge $v_1v_3$)
followed by a quadratic extension with base $v_1v_3$, a contradiction. 
Recall that $v_1$ is not adjacent to $v_3$ in $H$.

It remains to consider the case when both of the ends of $e$ have
degree three in $H$. Without loss of generality we assume that
$e=v_1v_2$, and hence $v_1,v_2,v_3$ and $v_4$ are pairwise distinct. It
follows that $G \setminus e$ is a strict linear extension of $H +
v_1v_3 + v_1v_4$ and is again a brick. This completes the case
analysis.
 ~\qed

\vskip 10pt

We conjecture the following strengthening of
Theorem~\ref{threevertices}.

\begin{conj} There exists $\alpha > 0$ such that every minimal brick $G$
has at least $\alpha |V(G)|$ vertices of degree three.
\end{conj}

Even a much weaker strengthening, namely, the conjecture that every
brick has at least four vertices of degree three, seems to require
new ideas or a substantial refinement of our techniques.

\junk{
\baselineskip 11pt
\vfill
\noindent
This material is based upon work supported by the National Science Foundation
under Grants No.~0200595 and~0354742. Any opinions, findings, and conclusions or
recommendations expressed in this material are those of the authors and do
not necessarily reflect the views of the National Science Foundation.
}

\end{document}